\documentclass[12pt]{article}
\usepackage[utf8]{inputenc}
\usepackage{amsmath, amssymb, amsthm}
\usepackage{xcolor}
\usepackage[colorlinks=true,
            linkcolor=black,
            urlcolor=blue!50!black,
            citecolor=blue!50!black]{hyperref}
\usepackage{hyperref}
\usepackage{geometry}
\begin{document}

\title{A New Primes-Generating Sequence}
\author{Mohammed Bouras}
\date{}

\maketitle

\begin{abstract}
For the sequence defined by
\[
a(n) = \frac{n^2 - n - 1}{\gcd\big(n^2 - n - 1,\, b(n-3) + n\,b(n-4)\big)}
\]
Where $b(n) = (n+2)\big(b(n-1) - b(n-2)\big),$ with initial conditions $b(-1) = 0$ and $b(0) = 1$, we find that $a(n)$ contains only $1$'s and primes, and can be represented as a finite continued fraction. It is more efficient for generating prime numbers than the Rowland sequence.
\end{abstract}


\section{Introduction}
In 2008, Rowland introduced an explicit sequence whose terms consist of $1$’s and prime numbers. This sequence is defined by the recurrence relation

\[
r(n) = r(n-1) + \gcd(n, r(n-1)), \quad r(1) = 7.
\] 
\\
Where $\gcd(x, y)$ denotes the greatest common divisor of $x$ and $y$.\\
The differences $r(n+1) - r(n)$ are \\
\\
1, 1, 1, 5, 3, 1, 1, 1, 1, 11, 3, 1, 1, 1, 1, 1, 1, 1, 1, 1, 1, 23, 3, 1, 1, \dots\ \text{(see \href{https://oeis.org/A132199}{A132199})} \\

While Rowland’s recurrence produces primes in the context of $\gcd$ additions, our sequence generates primes by filtering a quadratic polynomial through a rational function whose denominator is designed to divide out composite factors. The objective of this paper is to construct a new sequence that is more efficient for generating primes. We 
define the sequence using the $\gcd$ algorithm and the recurrence relation
\[
b(n) = (n+2)\big(b(n-1) - b(n-2)\big), \quad b(-1)=0, \quad b(0)=1.
\]
For all integers $n \geq 3$, the sequence is given by

\[
a(n) = \frac{n^2 - n - 1}{\gcd\big(n^2 - n - 1,\, b(n-3) + n\,b(n-4)\big)}.
\]
\\
Here, the numerator is a quadratic polynomial in $n$, while the denominator acts as a filtering mechanism that eliminates any composite factors shared with the sequence \( b(n-3) + n b(n-4) \). \\
\\
5, 11, 19, 29, 41, 11, 71, 89, 109, 131, 31, 181, 19, 239, 271, 61, 31, 379, 419, 461, 101, 29, 599, 59, 701, 151, 811, 79, 929, 991, 211, 59, 41, 1259, 1, 281, 1481, 1559, 149,
1721, 1, 61, 1979, 2069, 2161, 1, 2351, 79, 2549, 241, 1, 2861, 2969, 3079, 3191, 661, 311, 3539, 3659, 199, 71, 139, 4159, 4289, 4421, 911, 4691, 439, 4969, 269, 1051, 491,
179, 139, 5851, 1201, 101, 89, 1, 229, 1361, 6971, 1, 7309, 7481, 1531, 191, 8009, 431,
761, 1, 8741, 8929, 829, 9311, 1901, 109, 521, 10099, 10301, 191, 10711, 179, 359,
1031, 2311, 149, 631, 421, 401, 2531, 1171, 13109, 13339, 331, 251, 739, 131, 14519,
509, 3001, 151, 1409, 15749, 16001, 3251, 1, 409, 17029, 17291, 3511, 251, 18089,
1669, 601, 199, 19181, 1, 19739, 20021, 1, 349, 20879, 21169, 1951, 229, 22051, 22349,
1, 389, 4651, 23561, 23869, 24179, 1289, 1, 25121, 25439, 25759, 2371, 5281, 26731,
27059, 449, 1459, 181, 1, 28729, 709, 29411, 541, 971, 30449, 1621, 31151, 6301, 211,
1, 32579, 32941, 6661, 3061, 34039, \dots \ (\text{see \href{https://oeis.org/A356247}{A356247}})
\\

All terms are either equal to 1 or prime numbers, with no composite values observed among the first 10,000 terms verified computationally. Within this range, the value 1 appears exactly 1,420 times, accounting for approximately 14.2\% of the sequence, while the remaining 8,580 terms are primes.

It is immediate to observe that the combination \( b(n-3)+n b(n-4) \) can be replaced by \((n-1)!\) in the greatest common divisor without altering the result. The choice of the combination \( b(n-3)+n b(n-4) \) is preferable since it is typically much smaller than \((n-1)!\), which makes it more convenient for analysis and computation.

\section{Observations and Conjectures}
Let $x = n^2 - n - 1$ and $y = b(n-3) + n\,b(n-4)$. The behavior of the sequence 
$a(n) = \frac{x}{\gcd(x, y)}$ is determined by the common factors shared between \(x\) and \(y\). Three distinct cases occur:
\begin{itemize}
    \item \textbf{Coprime case} \(\gcd(x,y) = 1\): In this situation, the denominator shares no common factor with the quadratic numerator. Consequently, the sequence returns the full value of the quadratic
    $a(n) = n^2 - n - 1,$ which often yields a large prime.
    
    \item \textbf{Complete cancellation} \(\gcd(x,y) = x\): Here, the entire numerator is cancelled by the denominator, resulting in $a(n) = 1.$
    These are the only non-prime values in the sequence and occur precisely when the recursive expression \(y\) is a multiple of \(x\).
    
    \item \textbf{Partial cancellation} \(1 < \gcd(x,y) = d < x\): The sequence simplifies to $a(n) = \frac{x}{d},$
    which is still strictly greater than \(n\). Computational data confirm that these values are primes, with no known composite cases.
\end{itemize}
Together, these cases demonstrate that the sequence producing only 1's and primes. 
\\

We conjecture that

\begin{enumerate}
    \item Every term of this sequence is either $1$ or a prime number.
    \item The sequence contains all primes ending in $1$ or $9$.
    \item Except for $5$, the prime terms all appear exactly twice. Specifically, for any prime value $p = a(n)$, we also have

    \[
    a(p - n + 1) = p.
    \]
\end{enumerate}

Consequently, there exist integers \( n \) and \( m \) satisfying
\[
a(n) = a(m) = n + m - 1.
\]

Furthermore, we have
\[
a(n) = a(m) = \gcd(n^2 - n - 1,\, m^2 - m - 1).
\]
In this section, we generalize our result for the sequence \( a(n) \) derived from the finite continued fraction, as stated in the following theorem.

\section{Finite Continued Fraction Connection}
\textbf{Theorem 1.} Let $n \geq 3$ be an integer. Then the following identity holds:
\begin{equation}
\frac{mb(n-3) - nb(n-4)}{n(m-n+2) - m} \;=\;
\cfrac{1}{2 - \cfrac{3}{3 - \cfrac{4}{4 - \cfrac{5}{\ddots \; (n-1) - \frac{n}{m}}}}}
\tag{1}
\end{equation}

\textbf{Proof.} Consider the system of relations
\[
a_1 = 2a_2 - 3a_3, \quad a_2 = 3a_3 - 4a_4, \quad \ldots, \quad a_{n-1} = (n)a_n - (n+1)a_{n+1}.
\]

From the first equation, we obtain
\[
\frac{a_2}{a_1} = \frac{1}{\frac{2 a_2 - 3 a_3}{a_2}}.
\]

Proceeding recursively yields the finite continued fraction representation
\begin{equation}
\frac{a_{2}}{a_{1}} =
\cfrac{1}{2 - \cfrac{3}{3 - \cfrac{4}{4 - \cfrac{5}{\ddots \; (n-1) - \frac{n a_n}{a_{n-1}}}}}}
\tag{2}
\end{equation}

Comparing (1) and (2) immediately gives
\begin{equation}
m a_n = a_{n-1}.
\tag{3}
\end{equation}

Next, we express \( a_1 \) in terms of \( a_{n-1} \) and \( a_n \). Eliminating intermediate terms, we find
\begin{equation}
a_1 = (n-1)a_{n-1} - (n^2 - 2)a_n.
\tag{4}
\end{equation}

Substituting (3) into (4) yields
\[
a_1 = (n(m - n + 2) - m) a_n.
\]

Similarly, by tracing the recurrence for \( a_2 \), we obtain
\begin{equation}
a_2 = b(n-3) a_{n-1} - n b(n-4) a_n.
\tag{5}
\end{equation}

Using (3) in (5) gives
\[
a_2 = (m b(n-3) - n b(n-4)) a_n.
\]

Finally, substituting \( a_1 \) and \( a_2 \) into (2) exactly recovers the result in (1).

This completes the proof.
\\

\textbf{Theorem 2.} For all integers $n \geq 3$, the continued fraction
\[
\frac{2\left( m\,b(n-3) - n\,b(n-4) \right)}{n(m-n+1)}
=
\cfrac{1}{1 - \cfrac{1}{2 - \cfrac{2}{3 - \cfrac{3}{\ddots - \cfrac{n-1}{m}}}}}
\]

A special case of this identity relates directly to the left factorial function
$!n = \sum_{k=0}^{n-1} k!$ when $m = n$.

In this context, the auxiliary sequence $b(n)$ can be expressed as
\[
b(n) = (n+2)\,\frac{!(n+1)}{2}.
\]

\noindent\textbf{Proof.}  
Following the same reasoning as in Theorem 1, consider:
\[
a_1 = a_2 - a_3, \quad
a_2 = 2(a_3 - a_4), \quad
\ldots, \quad
a_{n-1} = (n-1)(a_n - a_{n+1}).
\]
By systematically substituting and simplifying these relations, we obtain the stated result.

\section{Main Results}
We define a family of sequences using a generalized denominator (Theorem 1).
\subsection{Quadratic Expression}
Consider the quadratic polynomial \( n^{2} + (k - 2)n - k \), and let \( m = -k \). The unreduced denominator of the finite continued fraction is as follows:

\[
a_k(n) = \frac{n^2 + (k-2)n - k}{\gcd\big(n^2 + (k-2)n - k,\, k\,b(n-3) + n\,b(n-4)\big)}.
\]

For sufficiently small \(n\), the sequence \(a_k(n)\), with \(k\) held fixed, predominantly produces large prime values arising from the associated quadratic form.

\begin{table}[h]
\centering
\caption{The sequence \( a_k(n) \) for \( k = 1, 2, 3, 4, 5 \).}
\label{tab:ak_sequence}
\begin{tabular}{|c|l|l|}
\hline
$k$ & $a_k(n)$ & OEIS \\
\hline
1 & 5, 11, 19, 29, 41, 11, 71, 89, 109, 131, 31, 181, 19, 239, 271, 61, 31, \dots & \href{https://oeis.org/A356247}{A356247} \\
2 & 7, 7, 23, 17, 47, 31, 79, 7, 17, 71, 167, 97, 223, 127, 41, 23, 359, 199, \dots & \href{https://oeis.org/A363102}{A363102} \\
3 & 3, 17, 9, 13, 53, 23, 29, 107, 43, 17, 179, 23, 79, 269, 101, 113, 29, 139, \dots & \href{https://oeis.org/A362086}{A362086} \\
4 & 11, 5, 31, 11, 59, 19, 19, 29, 139, 41, 191, 1, 251, 71, 29, 89, 79, 109, \dots & \href{https://oeis.org/A363347}{A363347} \\
5 & 13, 23, 7, 49, 13, 83, 103, 5, 149, 1, 29, 233, 53, 23, 67, 373, 59, 1, \dots & \href{https://oeis.org/A363482}{A363482} \\
\hline
\end{tabular}
\end{table}

A notable structural property of the sequence \( a_k(n) \) is its reflective symmetry. For every prime value \( p = a_k(n) \), there exists another index such that
\[
a_k(p - n - k + 2) = p,
\]
where \( k \) is a fixed integer. This duality holds for all observed terms and enables precise prediction of the positions of repeated primes. \\

For the sequence \( a_2(n) \), there exist prime numbers \( p \) such that \( p \) occurs exactly three times in the sequence. Moreover, if \( n \) denotes the index of the first occurrence of \( p \), then the index of the third occurrence satisfies
\[
a_2(p + n) = p.
\]

\subsection{Linear Combination}
We now turn to the linear form \((k+1)n - k\). Setting \(m = n + k\) yields the sequence
\[
a_k(n) = \frac{(k+1)n - k}{\gcd\big((k+1)n - k,\, b(n-2) + k\,b(n-3)\big)}.
\]

\begin{table}[ht]
\centering
\caption{The sequence \( a_k(n) \) for \( k = 1, 2, 3, 4, 5 \).}
\label{tab:ak_sequence}
\begin{tabular}{|c|l|l|}
\hline
$k$ & $a_k(n)$ & OEIS \\
\hline
1 & 5, 7, 3, 11, 13, 1, 17, 19, 1, 23, 1, 1, 29, 31, 1, 1, 37, 1, 41, 43, 1, 47, \dots & \href{https://oeis.org/A356360}{A356360} \\
2 & 7, 5, 13, 2, 19, 11, 5, 1, 31, 17, 37, 1, 43, 23, 1, 1, 1, 29, 61, 1, 67, 1, \dots & \href{https://oeis.org/A369797}{A369797} \\
3 & 3, 13, 17, 7, 5, 29, 11, 37, 41, 1, 7, 53, 19, 61, 1, 23, 73, 1, 1, 1, 89, 31, \dots & \href{https://oeis.org/A370726}{A370726} \\
4 & 11, 4, 7, 13, 31, 1, 41, 23, 17, 1, 61, 1, 71, 19, 1, 43, 1, 1, 101, 53, 37, \dots & \href{https://oeis.org/A372761}{A372761} \\
5 & 13, 19, 5, 31, 37, 43, 7, 11, 61, 67, 73, 79, 17, 1, 97, 103, 109, 23, 11, \dots & \href{https://oeis.org/A372763}{A372763} \\
\hline
\end{tabular}
\end{table}

For \(1 \leq k \leq 5\), all terms in the sequence \(a_k(n)\) consist of the integer 1 and the prime numbers, except for the unique occurrence of the value 4 in the sequence \(a_4(n)\).

\end{document}